\documentclass[color]{amsart}
\usepackage{fancybox,color,graphicx}

\usepackage[all,color,line]{xy}
\newtheorem{theorem}{Theorem}[section]
\newtheorem{corollary}[theorem]{Corollary}
\newtheorem{lemma}[theorem]{Lemma}

\theoremstyle{definition}
\newtheorem{definition}[theorem]{Definition}
\theoremstyle{remark}
\newtheorem{remark}[theorem]{Remark}
\newtheorem{example}{Example}[section]
\numberwithin{equation}{section}
\newcommand{\mx}{\mathbf x}
\newcommand{\mX}{\mathbf X}

\newcommand{\mA}{\mathbf A}

\newcommand{\mC}{\mathbf C}

\newcommand{\mW}{\mathbf W}

\newcommand{\la}{\lambda}
\newcommand{\tr}{{\rm tr}}\newcommand{\eps}{\varepsilon}
\newcommand{\bu}{\mathbf u}
\newcommand{\bx}{\mathbf x}
\newcommand{\bX}{\mathbf X}
\newcommand{\bY}{\mathbf Y}
\newcommand{\bZ}{\mathbf Z}
\newcommand{\bk}{\mathbf k}

\newcommand{\RR}{\mathbb{R}}
\newcommand{\CC}{\mathbb{C}}

\newcommand{\sC}{\mathcal{C}}
\newcommand{\Sn}{\mathcal{S}}
\newcommand{\WW}[1]{\mathcal{W}_N(#1)}
\usepackage{chicago}
\title[Polynomials in independent complex Wishart matrices]{Asymptotic normality for traces of
polynomials in independent complex Wishart matrices}

\author{W{\l}odzimierz  Bryc}
\thanks{\noindent Research partially supported by NSF
grant \#DMS-0504198.}
\address{
Department of Mathematical Sciences, University of Cincinnati, 2855
Campus Way, PO Box 210025, Cincinnati, OH 45221-0025, USA}
\email{Wlodzimierz.Bryc@UC.edu}

\date{Created: August 17, 2006; Revised: October 18, 2006; January 10, 2007; Printed: \today}
\begin{document}
%\centerline{\shadowbox{Preliminary {\bf confidential} draft. Not for citation, not for reference!}}

\maketitle

\begin{abstract}
We derive a
 non-asymptotic
expression for the
 moments of traces of monomials in several independent complex Wishart matrices,
 extending some explicit formulas  available in the literature.
We then deduce the explicit expression  for the cumulants.
From the latter, we read out the multivariate
normal approximation to the traces of finite families of polynomials in independent complex Wishart matrices.
\end{abstract}
\tolerance=1000
%\subection{Things to check}
%T. Kollo and D. von Rosen, Minimal moments and cumulants of symmetric matrices:
%an application to the Wishart distribution, J. Multivariate Anal. 55 (1995),149–164.
%\typeout{Things to check: Kollo-Rosen line 92}\typein{}
\section{Introduction}
%\subsection{Complex Wishart distribution}
The complex
Wishart distribution was introduced by Goodman \cite{Goodman-63}, and studied by numerous authors, starting with
 \cite{Goodman-63b,Khatri-65a,Khatri-65b,Srivastava-65}.
Following \cite{Graczyk-Letac-Massam-03}, we define the complex Wishart law by specifying its Laplace transform.
 \begin{definition} Let $\Sigma$ be an $N\times N$ Hermitian
 positive-definite matrix (i.e. $x^*\Sigma x>0$ for all non-zero
 $x\in\CC^N$), and $p\in\{1,2,\dots,N\}\cup(N,\infty)$.
 We will say that a random Hermitian  matrix $\mW$ has
 complex Wishart law with shape parameter $p$ and scale parameter $\Sigma$ if the joint Laplace transform of its entries
 is
 \begin{equation}
   \label{L}
   E\left(\exp(\tr(\theta \mW))\right)=(\det(I-\theta\Sigma))^{-p}
 \end{equation}
 for all Hermitian matrices $\theta$ such that $\Sigma^{-1}-\theta$ is positive-definite.

  \end{definition}
  We  write $\mW\in\WW{\Sigma,p}$. This notation is deliberately modeled on \cite[Section 7.2]{Anderson-03}, but
   the reader should keep in mind that
this paper is devoted solely to the complex case.

From Gindykin's theorem, see \cite{Graczyk-Letac-Massam-03}, we know that  \eqref{L}
 is indeed the Laplace transform of a probability measure which is  supported on the set of
 $N\times N$ Hermitian
 positive matrices, (i.e. on $\mA$ such that $x^*\mA x\geq 0$ for all $x\in\CC^N$). It is also well known
 that   $\mW\in\WW{\Sigma,p}$  with $p\in\mathbb{N}$ can be expressed
 as a quadratic form in a complex Gaussian matrix with i.~i.~d. entries, see \eqref{below}.

Asymptotic normality for traces of polynomials in a single real
Wishart matrix with $\Sigma=I$ is known since the work of
\cite{Arharov-71} and \cite[Theorem 4.1]{Jonsson-82}, who relied on combinatorial
methods similar to ours. More general unitary-invariant ensembles
have been studied starting with \cite{Diaconis-Shahshahani-94}; for
recent work and references, see
\cite{Collins-Mingo-Siniady-Speicher-06}. \cite{Bai-Silverstein-04}
proved asymptotic normality and derived elegant integral formulas for
the means and the variances of traces of analytic
 functions of a more general class of empirical covariances.
 Asymptotic normality for
multi-matrix %unitarily invariant
models
was studied in
\cite{Cabanal-Duvillard-01} who used stochastic analysis to
establish asymptotic normality for traces of polynomials in two
independent complex Wishart matrices with $\Sigma=I$, and   in
\cite{Mingo-Nica-04} and \cite{Kusalik-Mingo-Speicher-05} who
studied  polynomials in several independent
unitarily invariant Wishart matrices by combinatorial methods.

Our main result extends
\cite[Corollary 9.4]{Mingo-Nica-04} to more general covariances,
and is proved by the closely related combinatorial method.

\begin{theorem}\label{T3} Fix $\la_1,\dots,\la_s>0$ and  let $p_j(N)=\lceil\la_j N\rceil$, $j=1,\dots,s$.
Fix a sequence of positive-definite $N\times N$ matrices $\mC_1(N),\dots,\mC_s(N)$ such
that
 $$\lim_{N\to\infty}\frac1N\tr(q(\mC_1{(N)},\dots,\mC_s{(N)}))$$ exists for all non-commutative polynomials
  $q\in\mathbb{C}\langle\bx_1,\dots,\bx_s\rangle$.
Suppose
 $\mW_1^{(N)}\in\WW{\frac1N\mC_1{(N)},p_1(N)},\dots,
 \mW_s^{(N)}\in\WW{\frac1N\mC_s{(N)},p_s(N)}$ are independent.

  Given  polynomials
  $q_1, \dots, q_d\in\mathbb{C}\langle\bx_1,\dots,\bx_s\rangle$,
  define the $\mathbb{C}^d$-valued random variable
  $\bZ_N=(\xi_1(N),\dots,\xi_d(N))$ by
 $$\xi_j(N)=\tr(q_j(\mW_1^{(N)},\dots,\mW_s^{(N)})).$$
Then there is a  pair $(\bX,\bY)$ of  jointly Gaussian mean zero $\RR^d$-valued random variables  such that
 $$\bZ_N  -E(\bZ_N)\xrightarrow{\mathcal{D}} \bX+i\bY.$$
\end{theorem}
The statement of Theorem \ref{T3} does not specify the covariance matrix of the limiting normal law, and usually
the law of $\bZ=\bX+i\bY$ fails to be ``complex Gaussian", as the covariances of $\bX$ and $\bY$ may differ.
The expression for the joint $(2d)\times(2d)$ covariance of $(\bX,\bY)$ that can
be read out from the proof becomes more explicit in simple situations,
see Example \ref{Ex CLT} below.   More convenient representation of the covariance
could be of interest to statistics,
see \cite{Edelman-Mingo-RajRao-Speicher}.

We remark that if  $\la_j\geq 1$ for some $j$, we can take as shape parameter  $p_j=\la_j N$
instead of the first larger integer
$\lceil \la_j N\rceil$. We also note that under natural assumptions, the centering is linear in $N$,
see \eqref{lim mean}.

Theorem \ref{T3} is a consequence of  a number of non-asymptotic
results.
%The paper is organized as follows.
In Section \ref{Sect 1} we state the formula
 for mixed moments of traces of monomials in independent Wishart matrices, and we use it to derive
 some of the explicit expressions available in the literature.
In Section \ref{Sect Proof 1} we prove the formula for moments. In Section \ref{Sect Formal}, we  compute
 multivariate cumulants of the traces of monomials in independent Wishart matrices.
   In Section \ref{Sect CLT} we use the cumulants to prove Theorem \ref{T3}.
%Section \ref{Sect Asymp} \fbox{is not ready}.

\section{Mixed moments of independent complex Wishart matrices}\label{Sect 1}

Explicit formulas for moments of polynomials in one complex Wishart matrix appear in
\cite{Lu-Richards-01}, \cite{Hanlon-Stanley-Stembridge-92},
\cite{Graczyk-Letac-Massam-03}. Haagerup and Thorbj{\o}rnsen \cite{Haagerup-Thorbjornsen-03}
 compute the generating function
$E\tr(\exp s \mW)$ for Wishart matrices $\mW$ with scale parameter $\Sigma=I$.
Mingo and Nica \cite{Mingo-Nica-04} give explicit formulas for moments and cumulants
of monomials in several independent Wishart matrices of the same shape parameter $p$
and the same scale parameter $\Sigma=I$.
Here we adapt their technique to  more general shape and scale parameters.

\subsection{}%The setup}
We consider moments of monomials $q_{\pi,t,h}(\mx_1,\mx_2,\dots,\mx_s)$
of degree $n$
in non-commuting variables $\mx_1,\mx_2,\dots,\mx_s$, which are parameterized by
permutations $\pi\in\Sn_n$,  functions $t:\{1,\dots,n\}\to\{1,\dots,s\}$
and matrix-valued functions
$h:\{1,\dots,n\}\to\mathcal{M}_{n\times n}(\CC)$.
The interpretation of these parameters is as follows:
the cycles of permutation $\pi$  determine how the monomial is factored into the product of
traces; index $t(j)$
 describes which of the variables appears at the $j$-th position in the monomial,
 and $h(j)$ is an optional parameter which can be used, for example,
 to extract information about the moments of a single entry of the product.

 Let $\sC(\sigma)$ denote the set of cycles of permutation $\sigma$. Define
\begin{equation}
  \label{r}
q_{\sigma,t,h}(\mx_1,\mx_2,\dots,\mx_s)=\prod_{c\in \sC(\sigma)}\tr\left(\prod_{j\in c} h(j)\mx_{t(j)}\right).
\end{equation}
Here the notation $\tr\left(\prod_{j\in c} h(j)\mx_{t(j)}\right)$ is to be interpreted as
the trace of the product of the matrices taken in the same order in which the integers $j$ appear in the cycle $c$.
If $h\equiv I$ we will write $q_{\sigma,t,I}$; if $t\equiv 1$, we will write $q_{\sigma,1,h}$.

For example, %\eqref{Ex}
 $(\tr(\mW_1\mW_2))^2=q_{\sigma,t,I}(\mW_1,\mW_2)$ with
 \begin{equation}
   \label{Ex0}
   \sigma=(1,2),(3,4),\;  t(j)=(3+(-1)^j)/2.
 \end{equation}
Clearly, this representation is not unique; \cite{Hanlon-Stanley-Stembridge-92}
index their monomials
by the cycle type of permutation $\sigma$, i.e., by a partition $\la\vdash n$, see Corollary \ref{HLL}.
However, the statement of  Theorem \ref{T1}
relies on compositions of permutations, so the notation $q_{\pi,t,h}$ is more convenient here.
We also remark that monomials $q_{\sigma,t,h}$ are  natural generalizations to several variables of the
monomials $r_\sigma$ introduced
in \cite{Graczyk-Letac-Massam-03}, see
Corollary \ref{C-GLM}.

For a fixed $t:\{1,\dots,n\}\to\{1,\dots,s\}$, let
\begin{equation}
  \label{t-pres}
  \Sn_n(t)=\{\pi\in \Sn_n: t=t\circ\pi\}
\end{equation} be the set of
$t$-preserving permutations. Under the interpretation that $t$ assigns one of the $s$ colors
 to each of the integers $\{1,\dots,  n\}$,
these are "color-preserving"  permutations, and this is how we will refer to them.
Each $\alpha\in \Sn_n(t)$ may be identified with the $k$-tuple of permutations
of the $k$ nonempty sets among $t^{-1}(1), t^{-1}(2),\dots,t^{-1}(s)$. Accordingly, we  decompose
the set of cycles $\sC(\alpha)=\bigcup_{j=1}^s \sC_j(\alpha)$ into the union of sets $\sC_j(\alpha)$ of cycles on
$t^{-1}(j)$, $j=1,2,\dots,s$; here we allow the possibility that $t^{-1}(j)=\emptyset$ for some $j$, in which case
$\sC_j=\emptyset$.

If $\sigma,\alpha\in\Sn_n$, by $\sigma\alpha$ we denote the composition of permutations, i.e. the permutation
 $\sigma\circ\alpha$.

With the above notation, we  state the main result of this section.
\begin{theorem}
  \label{T1}
  If $\mW_1\in\WW{\Sigma_1,p_1},\mW_2\in\WW{\Sigma_2,p_2},\dots,\mW_s\in\WW{\Sigma_s,p_s}$ are independent  and $\sigma\in\Sn_n$, then
\begin{multline}
  \label{Main formula}
E\left(q_{\sigma,t,h}(\mW_1,\mW_2,\dots,\mW_s)\right)\\=\sum_{\alpha\in \Sn_n(t)}p_1^{\#\sC_1(\alpha)}p_2^{\#\sC_2(\alpha)}\dots
p_s^{\#\sC_s(\alpha)}q_{\sigma\alpha,t,h}(\Sigma_1,\Sigma_2,\dots,\Sigma_s).
\end{multline}
\end{theorem}
Theorem \ref{T1} follows by induction with respect to $s$ from the formula in
 \cite[Theorem 2]{Graczyk-Letac-Massam-03}, which we state as
Corollary \ref{C-GLM} below. However, the notation for such a proof is cumbersome, so
in Section \ref{Sect Proof 1} we give a direct elementary proof based on
 techniques of
\cite{Mingo-Nica-04}.

\subsection{}%Some special cases of interest}
We now use Theorem \ref{T1} to derive several explicit formulas which already have appeared
in the literature.
(The statements are re-cast into our notation.)

%\subsubsection{One matrix case}
\begin{corollary}[{\cite[Theorem 2]{Graczyk-Letac-Massam-03}}]\label{C-GLM}
For $\pi\in\Sn_n$, $h:\{1,\dots,n\}\to \mathcal{M}_{n\times n}(\CC)$, define
 $$r_{\pi}(h_1,\dots,h_n)(\mx)=\prod_{c\in \sC(\pi)}\tr\left(\prod_{j\in c}\mx h_j\right).$$
If $\mW\in\WW{\Sigma,p}$, then
  \begin{equation}
    E\left(r_{\pi_0}(h_1,\dots,h_n)(\mW)\right)=\sum_{\pi_1\in\Sn_n}p^{\#\sC(\pi_1^{-1}\pi_0)}r_{\pi_1}(h_1,\dots,h_n)(\Sigma).
  \end{equation}
\end{corollary}
\begin{proof}
Notice that $\pi_1$ and $\pi_1^{-1}$ have the same cycle type (i.e., their cycles have the same lengths), and
$$\tr(A_1A_2\dots A_k)=\overline{\tr(A_k^*\dots A_2^*A_1^*)}.$$
Since $t\equiv 1$, we have $\Sn_n(t)=\Sn_n$.
So with $\sigma=\pi_0^{-1}$, formula \eqref{Main formula} gives
\begin{equation*}
E\left( \overline{ r_{\pi_0}(h_1,\dots,h_n)(\mW)}\right)=E\left(q_{\sigma,1,h^*}(\mW)\right)
=\sum_{\pi\in\Sn_n}p^{\#\sC(\pi)}q_{\sigma\pi,1,h^*}(\Sigma).
\end{equation*}
Substituting $\pi_1^{-1}=\sigma\pi$ we get
\begin{multline*}
 E\left(\overline{ r_{\pi}(h_1,\dots,h_n)(\mW)}\right)=
\sum_{\pi_1\in\Sn_n}p^{\#\sC(\sigma^{-1}\pi_1^{-1})}q_{\pi_1^{-1},1,h^*}(\Sigma)\\
=\sum_{\pi_1\in\Sn_n}p^{\#\sC(\pi_0\pi_1^{-1})}\overline{r_{\pi_1}(h_1,\dots,h_n)(\Sigma)}.
\end{multline*}
This ends the proof after we recall that conjugation preserves the number of cycles:
since $\alpha\sigma$ and
 $\sigma\alpha$ are conjugated,
\begin{equation}
  \label{flip}
  \#\sC(\alpha\sigma)=\#\sC(\sigma\alpha).
\end{equation}
\end{proof}
\begin{corollary}[{\cite[Corollary 2.4]{Hanlon-Stanley-Stembridge-92}}]
  \label{HLL}
For $\la=(\la_1,\la_2,\dots)\vdash n$, let $p_\la(\mx)=\tr (\mx^{\la_1})\tr (\mx^{\la_2}) \dots$.
If $\mW\in\WW{I,p}$, and $\sigma\in\Sn_n$ has cycle type $\la$, then
\begin{equation}
  \label{HSS}
  E\left(p_\la(\mW)\right)=\sum_{\alpha\in \Sn_n} p^{\# \sC(\alpha)} N^{\# \sC(\alpha \sigma)}.
\end{equation}
\end{corollary}
\begin{proof} Since $p_\la(\mW)=q_{\sigma,1,I}(\mW)$ and $\tr(I)=N$, \eqref{HSS} follows
from \eqref{Main formula} by \eqref{flip}.
\end{proof}

%\subsubsection{Several matrices}
\begin{corollary}[{\cite[Lemma 9.1]{Mingo-Nica-04}}]
  \label{MN}
  If $\mW_1,\mW_2,\dots,\mW_s\in\WW{I,p}$ are independent, then
   \begin{equation}
     \label{MN-9.1}
     E(q_{\sigma,t,I}(\mW_1,\dots,\mW_s))=
     \sum_{\alpha\in \Sn_n(t)}p^{\#\sC(\alpha)}N^{\#\sC(\alpha^{-1}\sigma)}.
   \end{equation}
\end{corollary}
\begin{proof}
Using \eqref{flip},
$$ E(q_{\sigma,t,I}(\mW_1,\dots,\mW_s))= \sum_{\alpha\in\Sn_n(t)} p^{\#\sC(\alpha)} N^{\#\sC(\alpha\sigma)}=
\sum_{\alpha^{-1}\in\Sn_n(t)} p^{\#\sC(\alpha^{-1})} N^{\#\sC(\alpha^{-1}\sigma)},$$
which gives \eqref{MN-9.1}, as $\#\sC(\alpha)=\#\sC(\alpha^{-1})$.
\end{proof}

For asymptotic normality in Theorem \ref{T3}, it is interesting to determine
how the first moment of a trace depends on $N$ in a special setting.

 Fix $n\geq 1$, $\sigma=(1,2,\dots,n)\in\Sn_n$,
  $t:\{1,\dots,n\}\to\{1,\dots,s\}$, and $\la>0$.
  Let $\Sn_n^\circ(t)\subset\Sn_n(t)$ denote the set of color-preserving permutations $\alpha$ such that
  $\#\sC(\alpha)+\#\sC(\sigma\alpha)=n+1$.

\begin{corollary}\label{Col AM} Suppose $\mW_1^{(N)},\dots,\mW_s^{(N)}\in\WW{\mC_N/N,\lceil \la N\rceil }$ are
independent.  If there exist numbers $m_k$ such that
 \begin{equation}
   \label{C-conv}
   \lim_{N\to\infty}\tr(\mC_N^{k})-N m_k=0
 \end{equation}
 for all $k\in\mathbb{N}$, then
\begin{multline}
\label{lim mean}
\lim_{N\to\infty}\left(E\left(\tr(\mW_{t(1)}^{(N)}\dots \mW_{t(n)}^{(N)})\right)-
N\sum_{\alpha\in\Sn_n^\circ(t)} \la^{\# \sC(\alpha)}
\prod_{c\in \sC(\sigma\alpha)}m_{\# c}\right)=0.
%\\
%E\left(\tr(\mW_{t(1)}^{(N)}\dots \mW_{t(n)}^{(N)})\right)\\=
%N\sum_{\alpha\in\Sn_n^\circ(t)} \la^{\# \sC(\alpha)}
%\prod_{c\in \sC(\sigma\alpha)}m_{\# c}+O(1/N) \mbox{ as $N\to\infty$}.
\end{multline}
\end{corollary}
\begin{proof}
Denote $p(N)=\lceil \la N\rceil$. By Theorem \ref{T1},
\begin{multline*}
%\label{mean}
 E\left(\tr(\mW_{t(1)}^{(N)}\dots \mW_{t(n)}^{(N)})\right)\\
=N\sum_{\alpha\in\Sn_n(t)} \left(\frac{p(N)}{N}\right)^{\# \sC(\alpha)}N^{\#\sC(\alpha)+\#\sC(\sigma\alpha)-n-1}
\prod_{c\in \sC(\sigma\alpha)}\tr_N(\mC_N^{\# c}).
\end{multline*}
It is known, see  \cite{Jacques-68} and \cite[Theorem 1]{Machi-84}, that
\begin{equation}
  \label{genus}
\#\sC(\alpha)+\#\sC(\sigma)+\#\sC(\sigma\alpha)=n+2 - 2g
\end{equation} for some $g=g(\alpha)\in\{0,1,\dots,\lfloor(n-1)/2\rfloor\}$.
Clearly, $g(\alpha)=0$ for $\alpha\in \Sn_n^\circ(t)$.
Therefore,
\begin{multline*}
 E\left(\tr(\mW_{t(1)}^{(N)}\dots \mW_{t(n)}^{(N)})\right)
=N\sum_{\alpha\in\Sn_n(t)} \left(\frac{p(N)}{N}\right)^{\# \sC(\alpha)}N^{-2g(\alpha)}
\prod_{c\in \sC(\sigma\alpha)}\tr_N(\mC_N^{\# c})\\
=N\sum_{\alpha\in\Sn_n^\circ(t)} \left(\frac{p(N)}{N}\right)^{\# \sC(\alpha)}
\prod_{c\in \sC(\sigma\alpha)}m_{\# c}\\+
\sum_{\alpha\in\Sn_n^\circ(t)} \left(\frac{p(N)}{N}\right)^{\# \sC(\alpha)}
N\left(\prod_{c\in \sC(\sigma\alpha)}\tr_N(\mC_N^{\# c})-\prod_{c\in \sC(\sigma\alpha)}m_{\# c}\right)\\
+\sum_{g=1}^{\lfloor(n-1)/2\rfloor}\sum_{\alpha\in\Sn_n(t):g(\alpha)=g}
\left(\frac{p(N)}{N}\right)^{\# \sC(\alpha)}N^{1-2g}\prod_{c\in \sC(\sigma\alpha)}\tr_N(\mC_N^{\# c}).
\end{multline*}
Since $n$ is fixed, $p(N)/N\to \la$ and \eqref{C-conv} holds, this implies \eqref{lim mean}.

\end{proof}
\begin{remark} Inverses of permutations in $\Sn_n^\circ$
appear in \cite{Biane-97}, see the discussion of noncrossing permutations in \cite[Remark 2.10]{Mingo-Nica-04}.
\end{remark}
\begin{example}\label{Ex1} Under the assumptions of Theorem \ref{T1},
 formula \eqref{Main formula} used with $h(1)=E_{i,j}$, the matrix with exactly one non-zero element at $(i,j)$-th position,
 $h=I$ otherwise, gives
$$E\left(\mW_1\mW_2\dots\mW_s\right)=p_1p_2\dots p_s\Sigma_1\Sigma_2\dots\Sigma_s.$$
Indeed, we can identify the $(i,j)$-th entry of the matrix on the left hand side from
$E\left(\tr(E_{i,j}\mW_1\mW_2\dots\mW_s)\right)$. The only color-preserving permutation is the identity permutation.
Of course, the same answer arises by conditioning from the univariate case.
\end{example}
\subsection{Example of covariance calculation}
The purpose of this section is to  illustrate
 certain graphical representations which will be used in Theorem \ref{T2} and to derive formulas for Example \ref{Ex CLT}.

For complex-valued square-integrable random variables $\xi,\eta$, we denote
\begin{equation}
  \label{cov} \rm{cov}(\xi,\eta)=E(\xi\eta)-E(\xi)E(\eta);
\end{equation}
in this notation, $\rm{cov}(\xi,\xi)$ is a complex number.

{\tolerance=50000 Suppose $\mW_1,\mW_2$ are as in Theorem \ref{T1} so $\tr(\mW_1\mW_2)$  is real-valued. From Example \ref{Ex1} we know that
 $E(\tr(\mW_1\mW_2))=p_1p_2\tr(\Sigma_1\Sigma_2)$. To compute the variance of $\tr(\mW_1\mW_2)$ we therefore
 compute $E(\tr^2(\mW_1\mW_2))$ using \eqref{Main formula} with $n=4$ and $\sigma,t$ from \eqref{Ex0}.
 }
\begin{table}[h]
{\small
\begin{tabular}{|c||c|c|c|c|} \hline
$\alpha$
&(1)(2)(3)(4)&(1)(3)(2,4) &(1,3)(2)(4) &(1,3)(2,4)
\\
$\sigma\alpha$&(1,2)(3,4) &(1,2,3,4) & (1,4,3,2)&(1,4)(2,3)
\\
& & & &
\\
contribution&
 $p_1^2p_2^2\tr^2(\Sigma_1\Sigma_2)$ & $p_1^2p_2\tr((\Sigma_1\Sigma_2)^2)$&$p_1p_2^2\tr((\Sigma_1\Sigma_2)^2)$ &$p_1p_2\tr^2(\Sigma_1\Sigma_2)$
\\
\begin{minipage}{1.8cm}
\vspace{1cm} \centerline{ graph}
\end{minipage}  &\xymatrix  @-1pc{
{\color{blue}{\circ}\ar@{..}'[d]^{\color{blue}2}} &&\color{blue}\circ\ar@{..}[d]^{\color{blue}4} \\
{\bullet}\ar@{-}[d]_{1}& &{\bullet}\ar@{-}[d]_{3}  \\
\circ& &\circ
}
&
\xymatrix  @-1pc{
&{\color{blue}{\circ}\ar@{..}'[dl]^{\color{blue}2}} {\color{blue}\ar@{..}[dr]^4}& \\
{\bullet}\ar@{-}[d]_{1}& &{\bullet}\ar@{-}[d]_{3}  \\
\circ& &\circ
}
&
\xymatrix  @-1pc{
{\color{blue}{\circ}\ar@{..}'[d]^{\color{blue}2}} &&\color{blue}\circ\ar@{..}[d]^{\color{blue}4} \\
{\bullet}\ar@{-}[dr]_{1}& &{\bullet}\ar@{-}[dl]_{3}  \\
&\circ&
}
&
\xymatrix  @-1pc{
&{\color{blue}{\circ}\ar@{..}'[dl]^{\color{blue}2}} \ar@{..}[dr]^{\color{blue}4}& \\
{\bullet}\ar@{-}[dr]_{1}& &{\bullet}\ar@{-}[dl]_{3}  \\
&\circ&
}
\\
\hline
\end{tabular}\medskip
}
\caption{\label{Table1}}
\end{table}

Four possible $\alpha\in\Sn_4(t)$ and their contributions are listed in Table \ref{Table1}. The last row shows
 graphical representations
of the [airs $(\sigma,\alpha)$ which we will use in Section \ref{Sect Formal}: the cycles of $\sigma$ as cyclic clockwise
permutations of edges around the black vertices, and the cycles of $\alpha$ as cyclic clockwise
permutations of edges around the "colored" vertices marked as circles.
The first graph
 corresponds to
$(E(\tr(\mW_1\mW_2)))^2$, so the remaining connected graphs give the terms contributing to the variance.
Similar situation occurs for higher order cumulants, where again the contributions come only from the connected graphs.

Similarly, we can analyze $\tr(\mW_1\mW_2\mW_3)$. With $\sigma=(1,2,3)(4,5,6)$, there are
seven
terms contributing to  $\rm{cov}(\tr(\mW_1\mW_2\mW_3),\tr(\mW_1\mW_2\mW_3))$.
 Graphical representation of one such term with contribution $p_1p_2p_3\tr((\Sigma_1\Sigma_2\Sigma_3)^2)$
is depicted on the right in Fig. \ref{Fig1}.
We have therefore the following formulas.
\begin{example}[Example \ref{Ex1}, continued]\label{Ex2}If
$\mW_1\in\WW{\Sigma_1,p_1}$, $\mW_2\in\WW{\Sigma_2,p_2}$, $\mW_3\in\WW{\Sigma_3,p_3}$
are independent, then
  $$
  \rm{Var}\left(\tr(\mW_1\mW_2)\right)=p_1^2p_2^2\left(\frac{1}{p_1}+\frac{1}{p_2}\right)\tr((\Sigma_1\Sigma_2)^2)+p_1p_2\tr^2(\Sigma_1\Sigma_2);
  $$
  \begin{multline}\label{3W}
    %\rm{Var}\left(\tr(\mW_1\mW_2\mW_3)\right)=
    \rm{cov}(\tr(\mW_1\mW_2\mW_3),\tr(\mW_1\mW_2\mW_3))\\=
  p_1^2p_2^2p_3^2\left(\frac{1}{p_1}+\frac{1}{p_2}+\frac{1}{p_3}\right)\tr((\Sigma_1\Sigma_2\Sigma_3)^2)\\
  +p_1p_2p_3(p_1+p_2+p_3)\tr^2(\Sigma_1\Sigma_2\Sigma_3)+p_1p_2p_3\tr((\Sigma_1\Sigma_2\Sigma_3)^2).
  \end{multline}
\end{example}
%%%% calculation of cov
%\begin{tabular}{c|c|c}
%$\alpha$&$\sigma\alpha$ & contribution  \\ \hline
%(1, 4) (2, 5) (3, 6) & (6, 1) (4, 2) (5, 3) & $p_1p_2p_3\tr(\Sigma_1\Sigma_3)\tr(\Sigma_1\Sigma_2)\tr(\Sigma_2\Sigma_3)$ \\
%(1) (4) (2, 5) (3, 6) & (2, 4, 6, 1) (5, 3) & $p_1^2p_2p_3\tr(\Sigma_1\Sigma_2\Sigma_1\Sigma_3)\tr(\Sigma_2\Sigma_3)$\\
%(1, 4) (2) (5) (3, 6) & (6, 1) (3, 5, 4, 2) & $p_1p_2^2p_3\tr(\Sigma_3\Sigma_2\Sigma_1\Sigma_2)\tr(\Sigma_1\Sigma_3)$\\
%(1, 4) (2, 5) (3) (6) & (6, 5, 3, 1) (4, 2) & $p_1p_2p_3^2\tr(\Sigma_1\Sigma_3\Sigma_2\Sigma_3)\tr(\Sigma_1\Sigma_2)$\\
%(1) (4) (2, 5) (3) (6) & (2, 4, 6, 5, 3, 1) & $p_1^2p_2p_3^2\tr(\Sigma_1\Sigma_2\Sigma_1\Sigma_3\Sigma_2\Sigma_3)$\\
%(1) (4) (2) (5) (3, 6) & (2, 3, 5, 4, 6, 1) &$p_1^2p_2^2p_3\tr(\Sigma_1\Sigma_2\Sigma_3\Sigma_2\Sigma_1\Sigma_3)$\\
%(1, 4) (2) (5) (3) (6) & (6, 5, 4, 2, 3, 1) & $p_1p_2^2p_3^2\tr(\Sigma_1\Sigma_3\Sigma_2\Sigma_1\Sigma_2\Sigma_3)$
%\end{tabular}
%
%\medskip
In particular, if $p_1=p_2=p_3=p$ and $\Sigma_1=\Sigma_2=\Sigma_3=\Sigma$ then
$$  %\rm{Var}\left(\tr(\mW_1\mW_2\mW_3)\right)
\rm{cov}(\tr(\mW_1\mW_2\mW_3),\tr(\mW_1\mW_2\mW_3))=
  3p^5\tr(\Sigma^6)+
  +3p^4\tr^2(\Sigma^3)+p^3\tr(\Sigma^6).$$
A similar calculation yields
$$
\rm{cov}(\tr(\mW_1\mW_2\mW_3),\tr(\mW_3\mW_2\mW_1))=3p^5\tr(\Sigma^6)+3p^4\tr(\Sigma^2)\tr(\Sigma^4)+p^3\tr^3(\Sigma^2).
$$

%A similar calculation of $\rm{Var}(\tr(\mW^2)$ would involve $4!=24$ terms; graphical representation will allow us to select the
%terms that contribute to the limit as $N\to\infty$ under appropriate scaling of $p$ and $\Sigma$.
\section{Proof of Theorem \ref{T1}}\label{Sect Proof 1}

In the proof,  by %$A[i,j]$ or
$[A]_{i,j}$ we denote the $(i,j)$-th element of $A$.
It will also be convenient to write $h_j=h(j)$.

Following \cite{Mingo-Nica-04}, we write $q_{\pi,t,h}(\mx_1,\dots,\mx_s)$  as follows.
\begin{lemma}
  \label{L MingoNica}
\begin{equation}
  \label{MingoNica}
  q_{\pi,t,h}(\mx_1,\dots,\mx_s)=
\sum_{J:\{1,\dots,n\}\to\{1,\dots,N\}}
\prod_{i=1}^n [h_i \mx_{t(i)}]_{J(i),J(\pi(i))}.
\end{equation}
\end{lemma}
\begin{proof}
 For a single cycle $c=(i_1,i_2,\dots,i_w)$ of $\pi$,
%\begin{multline}\label{one cycle}
$$\tr\left(\prod_{r\in c}h_r\mx_{t(r)}\right)=
%\sum_{j_1,\dots,j_w=1}^N
%[h_1\mx_{t(1)}]_{j_1,j_2}[h_2\mx_{t(2)}]_{j_2,j_3}\dots [h_w\mx_{t(w)}]_{j_w,j_1}\\=
\sum_{J:c\to \{1,\dots,N\}}\prod_{i\in c}[h_{i}\mx_{t(i)}]_{J(i),J(\pi(i))}\;.$$
%\end{multline}
Thus
$$q_{\pi,t,h}(\mx_1,\dots,\mx_s)=
\sum_{J:\{1,\dots,n\}\to \{1,\dots,N\}}\prod_{c\in \sC(\pi)}\prod_{i\in c}[h_{i}\mx_{t(i)}]_{J(i)\;,
J(\pi(i))},
$$
and \eqref{MingoNica} follows.\end{proof}

\begin{proof}[Proof of Theorem \ref{T1}]
We first reduce the proof to the case when the shape parameters $p_1,\dots,p_s$ are positive integers. By independence and \eqref{L}, the joint Laplace transform
of $\mW_1,\dots,\mW_s$ is $$
  E\left(\exp(\sum_{r=1}^s\tr(\theta_r \mW_r))\right)=\prod_{r=1}^s\left(\det (I-\Sigma_r \theta_r)\right)^{-p_r}.$$
Denote $\theta_{r,i,j,0}=\Re[\theta_r]_{i,j}$, $\theta_{r,i,j,1}=\Im[\theta_r]_{i,j}$,
$W_{r,i,j,0}=\Re[\mW_r]_{i,j}$, $W_{r,i,j,1}=\Im[\mW_r]_{i,j}$.
Since
$$\tr(\theta_r\mW_r)=\sum_{i=1}^N \theta_{r,i,i,0}W_{r,i,i,0}+2\sum_{i=1}^N\sum_{j<i}(\theta_{r,i,j,0}W_{r,i,i,0}+\theta_{r,i,j,1}W_{r,i,j,1})$$
we therefore can compute the expectations of
 the entry-wise products
\begin{multline}\label{W-theta}
  \prod_{k=1}^n(2-\delta_{i_k,j_k})E\left(\prod_{k=1}^n W_{r_k,i_k,j_k,\eps_k}\right)\\=
\left.\frac{\partial^n}{\partial \theta_{r_1,i_1,j_1,\eps_1}\dots\partial \theta_{r_n,i_n,j_n,\eps_n}}
\prod_{r=1}^s\left(\det (I-\Sigma_r \theta_r)\right)^{-p_r}\right|_{\theta_1=\theta_2=\dots=\theta_s=0}
\end{multline}
for any $1\leq r_k\leq s$, $1\leq i_k,j_k\leq N$ and $\eps_k\in\{0,1\}$ such that $\eps_k=0$ when $i_k=j_k$.
As a function of $p_1,\dots,p_s$, the right hand side of \eqref{W-theta} is a polynomial.
 Since the left hand side of
\eqref{Main formula} is a linear combination of the expressions appearing on the left hand side of \eqref{W-theta},
 this shows that \eqref{Main formula} is an identity between polynomials so
it is enough to prove it for $p_1,p_2,\dots,p_s\in\mathbb{N}$.

We now prove \eqref{Main formula} with the assumption that  $p_1,p_2,\dots,p_s\in\mathbb{N}$. For $r=1,2,\dots,s$, let $\mX_r$ be a $p_r\times N$ matrix with i.~i.~d. complex  normal entries such that
 for  $i=1,2,\dots,p_r, j=1,2,\dots,N$ we have
\begin{equation}
  \label{complex moms}
  E\left([\mX_r]_{i,j}\right)=E\left([\mX_r]_{i,j}^2\right)=0,\; E\left(\left|[\mX_r]_{i,j}\right|^2\right)=1.
\end{equation} We also assume that matrices
$\mX_1,\dots,\mX_s$ are independent.
It is well known (\cite{Goodman-63}) that with $\Sigma_r=A_rA_r^*$, formula
\begin{equation}
  \label{below}
  \mW_r=A_r^*\mX_r^*\mX_rA_r
\end{equation}
defines independent
complex Wishart matrices  $\mW_r\in\WW{\Sigma_r,p_r}$, $r=1,\dots,s$.

Expanding matrix multiplications in \eqref{MingoNica}, we write the polynomial on the
left hand side of \eqref{Main formula} as
\begin{multline}\label{all cycles}
q_{\sigma,t,h}(\mW_1,\mW_2,\dots,\mW_s)\\
=\sum_{J,K,L,M}\left(\prod_{i=1}^n [\mX_{t(i)}^*]_{K(i),L(i)}[\mX_{t(i)}]_{L(i),M(i)}\right) \\
\times\prod_{i=1}^n
[h_iA_{t(i)}^*]_{J(i),K(i)}[A_{t(i)}]_{M(i),J(\sigma(i))}\;,
\end{multline}
where the sum is now taken over  all functions
$J,K,M:\{1,\dots,n\}\to\{1,\dots,N\}$ and over all  functions $L$ from the set
$$\mathcal{L}:=\{L:\{1,\dots,n\}\to\mathbb{N}: L(i)\leq p_{t(i)} \mbox{ for all $i$}\}.$$
Wick's formula expresses the expectation of the product of jointly normal centered
random variables as the sum  over all pairings of products of covariances. Under \eqref{complex moms},
the only contributing pairings are those that
 pair $[\mX_{t(i)}^*]_{K(i),L(i)}$ with $[\mX_{t(j)}]_{L(j),M(j)}$. Thus the
 contributing pairings are determined by the pairs $(i,\alpha(i))$ with
$\alpha\in\Sn_n(t)$,  and the expectation becomes
\begin{multline*}
E\left(\prod_{i=1}^n [\mX_{t(i)}^*]_{K(i),L(i)}[\mX_{t(i)}]_{L(i),M(i)}\right)\\=
\sum_{\alpha\in \Sn_n(t)}E\left(\prod_{i=1}^n \overline{[\mX_{t(i)}]}_{L(i),K(i)}
[\mX_{t(\alpha(i))}]_{L(\alpha(i)),M(\alpha(i))}\right)\\=
\sum_{\alpha\in \Sn_n(t)}1_{L=L\circ\alpha,M=K\circ\alpha^{-1}}.
\end{multline*}

This allows us to write the expectation of \eqref{all cycles} as
\begin{multline*}%\label{E(all cycles)}
E\left(q_{\sigma,t,h}(\mW_1,\mW_2,\dots,\mW_s)\right)\\
=\sum_{J,K,L,M}\sum_{\alpha\in \Sn_n(t)}1_{L=L\circ\alpha,M=K\circ\alpha^{-1}}
\prod_{i=1}^n
[h_iA_{t(i)}^*]_{J(i),K(i)}[A_{t(i)}]_{M(i),J(\sigma(i))}\\
=\sum_{J,K}\sum_{\alpha\in \Sn_n(t)}\sum_{L\in\mathcal{L}}1_{L=L\circ\alpha}
\prod_{i=1}^n
[h_iA_{t(i)}^*]_{J(i),K(i)}[A_{t(i)}]_{K(\alpha^{-1}(i)),J(\sigma(i))}
\end{multline*}
Since $t(j)=t(\alpha^{-1}(j))$ and $\sum_{L\in\mathcal{L}}1_{L=L\circ\alpha}=\prod_{r=1}^s p_{r}^{\#\sC_r(\alpha)}$,
this becomes \begin{multline*}
%\\=
\sum_{\alpha\in \Sn_n(t)}\prod_{r=1}^s p_{r}^{\#\sC_r(\alpha)}\sum_{J,K}
\left(\prod_{i=1}^n
[h_iA_{t(i)}^*]_{J(i),K(i)}\right)\prod_{j=1}^n [A_{t(\alpha^{-1}(j))}]_{K(\alpha^{-1}(j)),J(\sigma(j))}\\=
\sum_{\alpha\in \Sn_n(t)}\prod_{r=1}^s p_{r}^{\#\sC_r(\alpha)}\sum_{J,K}
\prod_{i=1}^n
[h_iA_{t(i)}^*]_{J(i),K(i)}[A_{t(i)}]_{K(i),J(\sigma(\alpha(i)))}
\\=
\sum_{\alpha\in \Sn_n(t)}\prod_{r=1}^s p_{r}^{\#\sC_r(\alpha)}\sum_{J}
\prod_{i=1}^n
[h_iA_{t(i)}^*A_{t(i)}]_{J(i),J(\sigma\alpha(i))}\\=
\sum_{\alpha\in \Sn_n(t)}\prod_{r=1}^s p_{r}^{\#\sC_r(\alpha)}\prod_{c\in \sC(\sigma\alpha)}
\tr\left(\prod_{i\in c}h_iA_{t(i)}^*A_{t(i)}\right)\\=
\sum_{\alpha\in \Sn_n(t)}\prod_{r=1}^s p_{r}^{\#\sC_r(\alpha)}\prod_{c\in \sC(\sigma\alpha)}
\tr\left(\prod_{i\in c}h_i\Sigma_{t(i)}\right)\\=
\sum_{\alpha\in \Sn_n(t)}\prod_{r=1}^s p_{r}^{\#\sC_r(\alpha)}q_{\sigma\alpha,t,h}(\Sigma_1,\Sigma_2,\dots,\Sigma_s).
\end{multline*}
\end{proof}

\section{Formal power series, cumulants, and graphical enumeration}\label{Sect Formal}
%In this section we consider a family of independent complex Wishart $N\times N$ matrices
% $\mW_1^{(N)},\dots,\mW_s^{(N)}$ with shape parameters
%$N \la_1>0,\dots,N\la_s>0$ and scale parameters $\frac1N\mC^{(N)}_1,\dots,\frac1N\mC^{(N)}_s$. To simplify notation,
%we will drop the superscript $N$ in the matrices; the normalized trace $\tr_N=\frac1N\tr$
% will occasionally serve as a reminder
%that all matrices depend on parameter $N$ through their dimension. Later on we will assume that
%$\lim_{N\to\infty}\tr_N(q(\mC_1^{(N)},\dots,\mC_s^{(N)}))$ exists for all polynomials $q$.

As a motivation, consider the integral
\begin{equation}
  \label{Z}
  Z_N(u_1,\dots,u_m):=E\left( e^{-NV_{\bu}(\mW_1,\dots,\mW_s)} \right),
\end{equation}
where
\begin{equation*}
  \label{V}
  V_\bu(\mW_1,\dots,\mW_s)=\sum_{j=0}^m u_j \tr(q_j(\mW_1,\dots,\mW_s)),
\end{equation*}
and $q_1,\dots,q_m\in\mathbb{C}\langle\bx_1,\dots,\bx_s\rangle$ are non-commutative monomials.
Of course, the integral might diverge, so it cannot be studied without
additional assumptions.
In view of recent progress \cite{Maurel-Segala-0608192} in multimatrix models,
it would be interesting to find the sufficient conditions for the asymptotic expansion as $N\to\infty$ for $\ln Z_N$.

In this section  our interest  is limited to the formal power series expansion for $\ln Z_N$,
which is well defined without any assumptions on the monomials, and
serves as a convenient notation to encode certain combinatorial relations.
The "formula" we give should not be interpreted as a statement about
 the actual function $Z_N$ that inspired the formal expansion. On the other hand,
 the formal power series  for $\ln Z_N$ is a well defined object which we will use to
prove asymptotic normality for the traces of polynomials.

 In order to avoid misinterpretation, we will place the tilde over all formal power series,
 even if they arise from  familiar analytic functions like $\ln(1+x)$, or $\exp$.

%\subsection{Background}
\subsection{Formal power series}
A ring of formal power series in $m$
variables  $u_1,\dots,u_m$ is the set of all
complex infinite sequences $\widetilde{a}=\{a_{k_1,\dots,k_m}: k_j\geq 0\}$, with the algebraic operations
  (linear combination, multiplication, differentiation, integration, and sometimes composition) that
  are defined as if they were inherited from
  the corresponding operations on the infinite power series
 $$\widetilde{a}(u_1,\dots,u_m)=\sum_{k_1,\dots,k_m\geq 0} u_1^{k_1}\dots u_m^{k_m} a_{k_1,k_2,\dots,k_m}.$$
 Well defined operations on formal power series are those that rely only on the finite number of the coefficients.
 Thus the $(k_1,k_2,\dots,k_m)$-th coefficient should not depend on the value of
 $M$ when the operations are performed on
 the polynomials of degree $M$ which were obtained by
 truncation $\mod u_1^M\dots u_m^M$ of a formal power series for all large enough $M\geq M_0(k_1,k_2,\dots,k_m)$.
 Such operations  act on infinite sequences $\widetilde{a},\widetilde{b}$
  regardless of the issues of convergence, and  the power series notation  serves only
  as a convenient mnemonic device.
 In particular, if $\widetilde{a}$ is a formal power series in one variable, then the
 composition $\widetilde{a}(\widetilde{b})$   is a well defined formal power series,
  whenever $b_{0,\dots,0}=0$.
  For example, if $\widetilde{a}(u)=\sum_{n\geq 1} n! u^n$, then
  $\widetilde{a}(\widetilde{a}(u))$ is well defined, and the first few terms are
  $\widetilde{a}(\widetilde{a}(u))=u+2\cdot 2! u^2+2( 3!+(2!)^2)u^3+\dots$.

To shorten the  notation, we will write
$$
\widetilde{a}(u_1,\dots,u_m)=\sum_{k_1,\dots,k_m\geq 0} \bu^{\bk} a_{\bk}.
$$
We will also write $\bk!=k_1!k_2!\dots k_m!$.

 Every infinitely-differentiable function $f(u_1,\dots,u_m)$ gives rise to the formal
 power series
 $$\widetilde{f}(u_1,\dots,u_m)=\sum_{k_1,\dots,k_m\geq 0}\frac{d_{\bk}}{\bk!}\bu^{\bk},$$
 where
 $$d_\bk=\left.\frac{\partial^{k_1}\dots\partial^{k_m}}{\partial u_1^{k_1}\dots \partial u_m^{k_m}}f\right|_{\bu=0}.$$
 Of course, the series does not have to converge,
 and even if it does,   it may fail to converge to $f(\bu)$.

\subsection{Cumulants}
The cumulants of real random variables $\xi_0,\xi_1,\dots,\xi_m$ with exponential
are often defined as the coefficients in the Taylor expansion
$$\ln E\exp(iu_1 \xi_1+i u_2\xi_2+\dots+i u_m\xi_m)=\sum_{k_1,\dots,k_m\geq 1}\frac{ (i\bu)^\bk}{\bk!}c_{\bk};$$
if $\xi_0,\xi_1,\dots,\xi_m$ have
 finite moments of all orders, then the left-hand side is differentiable at $0$, and the cumulants can be defined as the derivatives at $0$.
If a family of complex-valued random variables $\xi_0,\xi_1,\dots,\xi_m$
has all moments but the exponential moments fail to exist, the above approach will not work but
we can still define  the cumulants by their combinatorial relation to moments. This relation is encoded as
 an
identity between two formal power series.
Let $\widetilde{L}$ be the formal power series for $\ln(1+u)$. Let
$$\widetilde{M}(\bu)=\sum_{n=1}^\infty \frac{1}{n!}E((u_1 \xi_1+ u_2\xi_2+\dots+ u_m\xi_m)^n)=
\sum_{n\geq 1}\sum_{k_1+\dots+k_m=n} \frac{\bu^\bk}{\bk!}M_\bk,$$
where
$$M_\bk=E(\xi_1^{k_1}\dots\xi_m^{k_m}).$$
Then the composition $\widetilde{L}\left(\widetilde{M}(\bu)\right)$  is a well defined  formal power series;
formula
\begin{equation}\label{LM}
\widetilde{L}\left(\widetilde{M}(\bu)\right)=
\sum_{k_1+\dots+k_m\geq 1}\frac{ \bu^\bk}{\bk!}c_{\bk}
\end{equation}
encodes the relation between the moments on the left hand side and the cumulants on the right hand side.

We will need an explicit form of this relation  due to \cite{Leonov-Shirjaev-59}, see
\cite[Proposition 3.8.4]{Lando-Zvonkin-04}; we state in the form equivalent to
  \cite[(3.1)]{Speed-83}.
By $\mathcal{P}[1\dots n]$ we denote the set of all partitions $\mathcal{V}=\{B_1,\dots,B_v\}$
of $\{1\dots n\}$ into the disjoint non-empty sets $B_1,\dots,B_v$, which are called the blocks of $\mathcal{V}$.
Then the cumulants $c_{k_1,\dots,k_m}$ are determined recurrently in terms of the moment $M_\bk$
 and the lower order cumulants
as follows.
Let $n=k_1+\dots+k_m$ and let
$$(X_1,\dots,X_n)=(\begin{array}[t]{c}\underbrace{\xi_1,\dots,\xi_1}\\k_1\end{array},
\begin{array}[t]{c}\underbrace{\xi_2,\dots,\xi_2}\\k_2\end{array},\dots,\begin{array}[t]{c}\underbrace{\xi_m,\dots,\xi_m}\\k_m\end{array}).$$
To each nonempty $B\subset\{1,\dots,n\}$ we associate the corresponding cumulant
$$c_B:=c_{j_1,\dots,j_m},$$
with $j_r=\#\{u\in B: X_u=\xi_r\}$ for $1\leq r\leq m$. Of course, $c_\bk=c_{B}$ with $B=\{1,\dots,n\}$;
 other sets $B$ correspond to cumulants of order lower than $n$.

Then \eqref{LM} is equivalent to the following
\begin{equation}
  \label{Exp-Log}
  M_\bk=
\sum_{\mathcal{V}\in\mathcal{P}[1\dots n]}\prod_{B\in\mathcal{V}} c_B.
\end{equation}
%\subsubsection{Original version of Leonov-Shiryayev formula} With $M_\bk=E(\xi_1^{k_1}\dots\xi_m^{k_m})$,
%\begin{equation}
%  \label{leonov-Shiryaev}
%  \frac{M_\bk}{\bk!}=\sum_{s\geq 1}\frac1{s!}\sum_{\bj_1+\dots+\bj_s=
%  \bk, \bj_i\ne0}\prod_{i=1}^s\frac{c_{\bj_i}}{\bj_i!}.
%\end{equation}
For example, with $m=2$ we have $c_{1,0}=E(\xi_1)$, $c_{0,1}=E(\xi_2)$, and
$c_{1,1}=\rm{cov}(\xi_1,\xi_2)$, consistently with notation in \eqref{cov}.

Since the cumulants determine moments, and the moments of real random variables determine the
distribution, provided they do not grow too fast, we get the following.
\begin{lemma} \label{Marcinkiewicz}
If $\xi_1,\dots,\xi_m$ are real-valued random variables with cumulants such that
$c_{\bk}=0$ for all $k_1+\dots+k_m\geq 3$ ,
then the joint law of $(\xi_1,\dots,\xi_m)$ is Gaussian.
\end{lemma}
We remark that by \cite{Marcinkiewicz}, condition $k_1+\dots+k_m\geq 3$
can be replaced by $k_1+\dots+k_m\geq M$ for some fixed $M$. We also note that Lemma \ref{Marcinkiewicz}
does not hold for complex random variables:
if $\xi,\eta$ are independent, $\xi$ is $N(0,1)$, and $\eta$ has uniform distribution
on $|z|=1$, then $\xi+\eta$ has the same moments, and hence the same
cumulants, as the $N(0,1)$ law.

\subsection{Maps and genus} Our expression for the cumulants uses sums over the hypermaps. We need
to recall background information and introduce the notation. We will follow \cite{Lando-Zvonkin-04} when possible; however
we need to clarify the role of coloring, and our formulas use the opposite orientation on the faces.

\subsubsection{Motivation: maps and pairs of permutations}
A topological (oriented) map $M$ on a compact oriented surface $S$ without boundary is a partition of $S$ in
three finite sets of cells:\begin{itemize}
\item[(i)] The set $V=V(M)$ of the vertices of $M$, which is a finite set of points;
\item[(ii)] The set $E=E(M)$ of the edges of $M$, which is a finite set of simple opened Jordan arcs, disjoint in pairs, whose
extremities are vertices;
\item[(iii)] The set $F=F(M)$ of the faces of $M$. Each face is homeomorphic to an open disc,
and its boundary is a union of vertices
and edges.
\item[(iv)] The rooted  map has a selected oriented edge, i.e., an edge with
 an adjacent vertex.
\end{itemize}
Since the graph $(V,E)$ is embedded into $S$,
the orientation of $S$ determines the unique cyclic order on the edges adjacent to each vertex.
(In the drawings, we will always follow the clockwise direction.) We note that the same graph can be embedded into
$S$ in many ways, see \cite[Fig. 1.17]{Lando-Zvonkin-04}.

The genus $g=0,1,2,\dots$ of a map $M$ is the genus of surface $S$, and is determined from Euler's formula
\begin{equation}
  \label{Euler}
  \#V-\#E+\#F=2-2g.
\end{equation}
(In particular, all maps on $S$ have the same genus.)

A rooted combinatorial map is an ordered pair of permutations $(\sigma,\alpha)\in \Sn_n\times\Sn_n$
acting transitively on $\{1,\dots,n\}$  such that $\alpha$ is an involution with no fixed points.
Such a pair is identified with the topological map $M=(V,E,F)$ through the associated sets of cycles,
with
$V=\sC(\sigma)$, $E=\sC(\alpha)$, and $F=\sC(\sigma\alpha)$.
As the rooted edge we take the cycle  $(1,\alpha(1))$; as distinguished vertex
we take the cycle of $\sigma$ which contains $1$. It is known that $(V,E,F)$ can be embedded into an oriented
 surface $S$, see
\cite[page 35]{Lando-Zvonkin-04}.

Conversely, to any a graph $(V,E)$ with fixed cyclic order on the edges adjacent to each vertex,
we can assign a pair of permutations as follows: we color all vertices "black",
split each edge in two by inserting a "white" vertex in the middle. We then label all new edges by numbers
$1,\dots,n$ with $n=2\#V$. (For rooted maps, we label with $1$ the distinguished (half)
edge adjacent to the distinguished vertex.)
Define $\sigma\in\Sn_n$ as a permutation whose cycles are labels of edges adjacent to each black vertex
taken in the cyclic order inherited from the map. Define $\alpha$ as the permutation
whose cycles are labels of edges adjacent to each white vertex.

We remark that the standard notation for the permutations $(\sigma,\alpha)$ of a map
is $(\sigma,\alpha,\phi)$ with $\phi=\alpha^{-1}\sigma^{-1}$ representing the faces;
our use of $\sigma\alpha$ instead of $\phi$
yields the opposite orientation on the faces and is more convenient in formulas below.

\subsubsection{Hypermaps with coloring}
\begin{definition}
A hypermap $H$ is any ordered pair of permutations $(\sigma,\alpha)\in \Sn_n\times\Sn_n$ which acts transitively on $\{1,\dots,n\}$.
 The genus $g=g(H)$   is defined by \eqref{genus}.
%\begin{equation}
%  \label{genus}
%\#\sC(\alpha)+\#\sC(\sigma)+\#\sC(\sigma\alpha)-n=2-2g.
%\end{equation}
\end{definition}
Clearly, $2g\leq n-1$, and by a result of \cite{Jacques-68},   $g$ is a non-negative integer.
(For an  accessible algebraic proof see \cite[Theorem 1]{Machi-84}.
For a graph-theoretic proof, see \cite[Proposition 1.5.3]{Lando-Zvonkin-04}.)

We will actually use colored hypermaps, where $\alpha\in\Sn_n(t)$. Following  \cite[Section 1.5]{Lando-Zvonkin-04}, we will identify each hypermap with a bipartite rooted colored
map
$H=(V=V_0\cup \bigcup_jV_j,E,F)$ whose
black  vertices are
$V_0=\sC(\sigma)$, and colored vertices of color $j$ are $V_j= \sC_j(\alpha)$.
The  edges are $E=\{1,\dots,n\}$ and the faces are $F=\sC(\sigma\alpha)$.

This identification can also be described as follows. Given a colored connected bipartite map
$H=(V=V_0\cup \bigcup_jV_j,E,F)$ in which
 colored edges connect vertices in $V_0$ to the colored vertices of matching color
  in $\bigcup_{j\geq 1}V_j$,
 we identify (label) $E=\{1,\dots,n\}$; the root edge is labeled $1$. We define
$\sigma$ as a cyclic
permutation, in clockwise direction,  of the edges attached to each vertex  in $V_0$.
We define permutation $\alpha$  as the
cyclic clockwise
permutation of the edges attached to
each of the colored vertices. Thus $\sigma,\alpha\in\Sn_E$, and $\alpha$ is color-preserving.
Under this identification, formula \eqref{genus} is just Euler's formula, so $g\in\{0,1,2,\dots,\lfloor(n-1)/2\rfloor\}$.

We now introduce the bipartite maps
 with stars of type $q_1,q_2,\dots q_m$; this terminology is motivated by
 \cite{Guionnet-Maurel-Segala-0503064}.
 \begin{definition}
 Let $q_1,q_2,\dots,q_m$ be (ordered) sequence of monomials of degrees $n_1,\dots,n_m\geq 1$
 in non-commutative variables $\bx_1,\dots,\bx_s$. Let
$n=n_1+\dots+n_m$ denote the degree of the product
$q_1q_2\dots q_m$. Let  $t:\{1,\dots,n\}\to\{1,\dots,s\}$ be the coloring of $\{1,\dots,n\}$
  according to the labels
 of variables  $\bx_1,\dots,\bx_s$ in the monomial $q_1\dots q_m$.
The set of
 hypermaps $\mathcal{H}(q_1,q_2,\dots,q_m)$ with stars of type $(q_1,q_2,\dots q_m)$,
 consist of pairs $(\sigma,\alpha)$
 such that\begin{enumerate}
 \item
 permutation $\sigma\in\Sn_n$ is fixed as
 %, is fixed with $m$ cycles of lengths corresponding to the degrees of
% $q_1,\dots,q_m$ by
 $$
 \sigma=(1,2,\dots,n_1)(n_1+1,\dots,n_1+n_2),\dots,(n_1+\dots+n_{m-1}+1,\dots,n_1+\dots+n_{m}).
 $$
 \item  permutations $\alpha\in\Sn_n(t)$
vary among the color-preserving permutations such that
 the action of the group $\langle\sigma,\alpha\rangle$ on $\{1,\dots,n\}$ is transitive.
\end{enumerate}
\end{definition}

In graph language, $\mathcal{H}(q_1,q_2,\dots,q_m)$ is the set of rooted bipartite colored maps $(V=V_0\cup V_c,E,F)$ with
$m$ black vertices $V_0$ of degrees $\deg q_1,\dots,\deg q_m$ corresponding to the monomials $q_1,q_2,\dots q_m$.
The
$n$ edges are labeled, cyclically ordered at each black vertex, and colored with different colors
associated to each variable $\bx_i, i\in\{1,\dots,s\}$, colored clockwise in the same order in which they appear in the monomials. The
number of the colored vertices $V_c$ and the ordering of their edges are not specified and vary subject to the following constraints:
 every edge connects a black vertex to  a colored vertex, and
only edges of the same color may connect to the colored vertex. The root edge corresponds to the first variable in $q_1$.

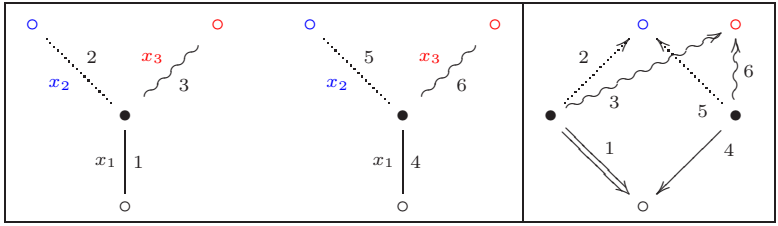
\begin{figure}[hbt]
\begin{tabular}{|c|c|}
\hline \xymatrix{
%{\color{blue}{\circ}\ar@{-}'[d]_{\color{blue}x_2}^2} &&\color{blue}\circ\ar@{-}[d]_{\color{blue}x_2}^4 &
%\color{red}\circ\ar@{-}[d]_{\color{red}x_3}^6
\color{blue}{\circ}\ar@{..}'[dr]_{\color{blue}x_2}^2&&\color{red}{\circ}\ar@{~}'[dl]_{\color{red}x_3}^3
&\color{blue}{\circ}\ar@{..}'[dr]_{\color{blue}x_2}^5&&\color{red}{\circ}\ar@{~}'[dl]_{\color{red}x_3}^6\\
&{\bullet}\ar@{-}[d]_{x_1}^1& & &{\bullet}\ar@{-}[d]_{x_1}^4
\\
&\circ& & &\circ
}
&
\xymatrix{
&{\color{blue}{\circ}\ar@{<..}[dl]_{2}}\ar@{<..}[dr]_>>{5} &\color{red}\circ\ar@{<~}[dll]^>>>>>>>{3}\\
{\bullet}\ar@2{->}[dr]^{1}& &{\bullet}\ar@{->}[dl]^<<4{} \ar@{~>}[u]_6{}  \\
& \circ&
}
\\ \hline
\end{tabular}
\caption{On the left: two stars of type $\bx_1\color{blue}\bx_2\color{red}\bx_3$ for calculation of \eqref{3W}. Here
 $\sigma=(1,2,3)(4,5,6)$
On the right: the rooted colored hypermap of genus 1  with
 $\alpha=(1, 4) (2, 5) (3, 6)$, $\sigma\alpha=(1,5, 3, 4, 2, 6)$.
The edges are oriented from the black to the colored vertices, and the root edge is marked by the double arrow.
This hypermap has one face with the boundary
${\bf 1}\underline{4}{\bf 5}\underline{2}{\bf 3}\underline{6}{\bf 4}\underline{1}{\bf 2}\underline{5}{\bf 6}\underline{3}$;
 the underlined edges have negative orientation and do not appear in the cycle representation.\label{Fig1}}

\end{figure}

All edges in the graph have natural orientation from the black vertex to the colored vertex.
The cycles of permutation
$\sigma$ list edges adjacent to a given black vertex in clockwise direction. The cycles of
permutation $\alpha$ give the clockwise order of
the edges adjacent to a colored vertex; all these edges are required to have the same color. Then
the cycles of $\sigma\alpha$ list those edges of a face that are in counter-clockwise (positive)
orientation with respect to that face,
see Fig \ref{Fig1}.

 By $V_j$ we denote the set of all vertices in $V_c$ linked to an edge of color
 $j\in\{1,\dots,s\}$. We will also need the following special notation.
\begin{definition}
  \label{Def-notation}
  If $f$ is a face corresponding to a cycle $c=(j_1,j_2,\dots,j_k)$ of $\sigma\alpha$, by
  $\tr_N\mC^{f,\,t}$ or by $\tr_N\mC^{c,\,t}$ we denote
  $\frac1N\tr\left(\mC_{t(j_1)}\mC_{t(j_2)}\dots \mC_{t(j_k)}\right)$.
\end{definition}

\subsection{Genus expansion for the partition function}

\begin{definition}We define the formal power series  $\widetilde{\ln Z}_N(\bu)$
 as the composition of the formal power series $\widetilde{L}(\widetilde{Z}_N(\bu)-1)$,
where  $\widetilde{L}(z)=\sum_{n=1}^\infty (-1)^{n+1} \frac{z^n}{n}$ is the formal power series
of the function $L(z)=\ln(1+z)$, and
\begin{multline}
  \label{Z -formal}
  \widetilde{Z}_N(\bu)=\sum_{k_1,\dots,k_m\geq 0} N^{k_1+\dots+k_m}\frac{(-\bu)^{\bk}}{\bk!}
E\left(\prod_{j=1}^m\tr^{k_j}(q_j(\mW_1,\dots,\mW_s))\right).
\end{multline}
\end{definition}
Of course, the latter series is the formal power series of directional derivatives at zero
of function $Z_N(\bu)$, when
the monomials $q_1,\dots q_m$
are such that the integral \eqref{Z} converges for $u_1,\dots,u_m>0$.
Even if the series fails to converge, the coefficients of $\widetilde{\ln Z}_N(u_1,\dots,u_m)$ are
the multivariate cumulants of the family of random variables
$$\xi_1=N\tr (q_1(\mW_1,\dots,\mW_s)),\dots,\xi_m=N\tr (q_m(\mW_1,\dots,\mW_s)).$$

Let
$\mathcal{H}(g,\bk)$  denote the set of hypermaps of genus $g$ with $k_j$ stars of type $q_j$, $j=1,\dots,m$. Note that since $\bk$ is an ordered $m$-tuple,
all hypermaps in $\mathcal{H}(0,\bk)\cup\mathcal{H}(1,\bk)\cup\dots$ share the same $\sigma=\sigma_\bk$ and coloring $t=t_\bk$, and
they differ only in the choices of $\alpha$.
\begin{theorem}\label{T2} If $\mW_1\in\WW{\frac1N\mC_1,N\la_1},\dots,\mW_s\in\WW{\frac1N\mC_s,N\la_s}$ are independent, then
\begin{equation}
  \label{log Z}
  \frac{1}{N^2}\widetilde{\ln Z}_N(\bu)=\sum_{g\geq 0}\frac{1}{N^{2g}}
  \sum_{k_1+\dots+k_m\geq 1}\frac{(-\bu)^{\bk}}{\bk!}c(\bk,g,\mC),
\end{equation}
where
$$c(\bk,g,\mC)=
\sum_{(V,E,F)\in\mathcal{H}(g,\bk)}\prod_{r=1}^s \la_r^{\# V_r} \prod_{f\in F}\tr_N\mC^{f,\,t}.$$

\end{theorem}
\begin{proof}
The idea of proof is fairly simple to describe. By Theorem \ref{T1}, formula \eqref{Z -formal}
can be reinterpreted as a generating function of disjoint unions of hypermaps.
It is well known in combinatorics that a logarithm of such a generating function is the generating function
 for the connected objects of the same type, provided that the coefficients are multiplicative.
  More precise statement relies on cumbersome notation,
see the last two paragraphs that precede \cite[Proposition 3.8.3]{Lando-Zvonkin-04};
a simpler version of the multiplicative property in
\cite[Corollary 5.1.6]{Stanley-99v2} is too strong for our needs, so we will verify directly
 \eqref{Exp-Log}.

For $k_1+\dots+k_m\geq 1$, denote
$$
M_\bk=E\left(\prod_{j=1}^m\tr^{k_j}(q_j(\mW_1,\dots,\mW_s))\right).
$$
With $n_j=\deg(q_j)\geq 1$, denote by  $n_\bk=k_1n_1+\dots+k_mn_m$ the degree of the
monomial $q_1^{k_1}\dots q_m^{k_m}$, and let
 $t_\bk:\{1,\dots,n_\bk\}\to\{1,\dots,s\}$  be the mapping identifying the variables in this monomial.
%Fix a permutation $\sigma\in \Sn_{n_\bk}$ with $k_j$ cycles of length $n_j$, $j=1,\dots,m$.
Let $\sigma\in \Sn_{n_\bk}$ be the permutation with $k_j$ cycles of length $n_j$, $j=1,\dots,m$, so that
$\prod_{j=1}^m\tr^{k_j}(q_j(\bx_1,\dots,\bx_s))=q_{\sigma,t_\bk,I}(\mx_1,\dots,\mx_s)$ as in \eqref{r}.
%$$(1,\dots,n_1)(n_1+1,\dots,2n_1),\dots,((k_1-1)n_1+1,\dots,n_1k_1),(n_1k_1+1,\dots,n_1k_1+n_2),\dots,(
%of type $\langle n_1^{k_1}\dots n_m^{k_m}\rangle$.
From \eqref{Main formula} we have
\begin{multline*}
  M_\bk=
  E\left(q_{\sigma,t_\bk,I}(\mW_1,\dots,\mW_s)\right)\\=
  \sum_{\alpha\in\Sn_{n_\bk}(t_\bk)}\prod_{r=1}^s(\la_r N)^{\#\sC_r(\alpha)}N^{\#\sC(\sigma\alpha)-n_\bk}\prod_{c\in \sC(\sigma\alpha)}\tr_N \mC^{c,\,t}
\\=
  \sum_{\alpha\in\Sn_{n_\bk}(t_\bk)}N^{\#\sC(\alpha)+\#\sC(\sigma\alpha)-n_\bk}\prod_{r=1}^s \la_r ^{\#\sC_r(\alpha)}
  \prod_{c\in \sC(\sigma\alpha)}\tr_N \mC^{c,\,t}.
\end{multline*}
Combining this with \eqref{Z -formal} and using the fact that $\#\sC(\sigma)=k_1+\dots+k_m$, we get
\begin{equation*}
 M_\bk= \sum_{\sigma\in \Sn_{n_\bk}}\sum_{\alpha\in \Sn_{n_\bk}(t_\bk)}N^{\chi(\sigma,\alpha)}
  \prod_{r=1}^s \la_r ^{\#\sC_r(\alpha)}\prod_{c\in \sC(\sigma\alpha)}\tr_N \mC^{c,\,t},
\end{equation*}
where $\chi(\sigma,\alpha)=\#\sC(\sigma)+
  \#\sC(\alpha)+\#\sC(\sigma\alpha)-n_\bk$.
We note that by \eqref{genus}, we have $\chi(\sigma,\alpha)=2-2g$ when the group $\langle\sigma,\alpha\rangle$
generated by $\sigma,\alpha$ acts transitively on $\{1,\dots,n_\bk\}$.
  In general, we decompose $\{1,\dots, n_\bk\}$ into the orbits of $\langle\sigma,\alpha\rangle$.
 Since the edges of each star of type $q_j$ are in the same orbit of $\langle\sigma,\alpha\rangle$, there are at most $n:=k_1+\dots+k_m$
   orbits, which
   are determined by the following procedure.
   \begin{itemize}
   \item Fix the sequence of polynomials (stars)
   \begin{equation}
\label{QQQ}
   (Q_1,Q_2,\dots,Q_n)=(\begin{array}[t]{c}\underbrace{q_1,\dots,q_1}\\k_1\end{array},
\begin{array}[t]{c}\underbrace{q_2,\dots,q_2}\\k_2\end{array},\dots,\begin{array}[t]{c}\underbrace{q_m,\dots,q_m}\\k_m\end{array}).
   \end{equation}
   \item Choose a partition $\mathcal{V}=\{B_1,\dots,B_v\}$   of $\{1,\dots,n\}$. The blocks  of $\mathcal{V}$ represent
   the  clusters of connected stars with  no connections between  different blocks.
    Each such choice  splits $\sigma\in\Sn_{n_\bk}$ into
    the uniquely determined permutations $\sigma_B$ of the edges of stars assigned to block
     $B\in\mathcal{V}$.
\item Sum over all color-preserving permutations $\alpha\in\Sn_{n_\bk}(t)$
which connect the stars within each block at least once. Each such permutation can be split into the sequence of
color-preserving permutations $\alpha_B$ of
the edges of stars in block $B\in\mathcal{V}$, and $\langle\sigma_B,\alpha_B\rangle$ acts transitively on
the set of edges of  the stars in $B$.
  \end{itemize}
For $B\subset\{1,\dots,n\}$, we denote by $\mathcal{H}(B)$ the set of all hypermaps
  with stars of type $\{Q_j : j\in B\}$, which we order according to their position in \eqref{QQQ}.
  Let $\Sn_{B}^T(t)$ be the color preserving permutations of
the edges of stars in block $B\in\mathcal{V}$ such that $\langle\sigma_B,\alpha_B\rangle$ acts transitively.
  Using this notation, we have
\begin{multline*}
  M_\bk=\sum_{\mathcal{V}\in\mathcal{P}[1\dots n]}\prod_{B\in\mathcal{V}}
  \sum_{\alpha\in \Sn_{B}^T(t)}N^{\chi(\sigma,\alpha)}
  \prod_{r=1}^s \la_r ^{\#\sC_r(\alpha)}\prod_{c\in \sC(\sigma_B\alpha)}\tr_N \mC^{c,\,t}\\
  \\=\sum_{\mathcal{V}\in\mathcal{P}[1\dots n]}\prod_{B\in\mathcal{V}}
   \sum_{H=(V,E,F)\in\mathcal{H}(B)}N^{2-g(H)}\prod_{r=1}^s\la_r^{\#V_r}\prod_{f\in F}\tr_N \mC^{f,\,t}.
  \end{multline*}
Comparing this to \eqref{Exp-Log} with $\xi_j=\tr(q_j(\mW_1,\dots,\mW_s))$, $j=1,\dots,m$,
 we see that the corresponding cumulants are
$$c_\bk= \sum_{g=0}^{\lfloor(n-1)/2\rfloor}\sum_{H=(V,E,F)\in\mathcal{H}(g,\bk)}N^{2-g}\prod_{r=1}^s\la_r^{\#V_r}
\prod_{f\in F}\tr_N \mC^{f,\,t},$$
so
\begin{multline*}
  \frac{1}{N^2}\widetilde{\ln Z}_N(\bu)=\sum_{n\geq 1}\sum_{k_1+\dots+k_m=n}\frac{(-\bu)^\bk}{\bk!}c_\bk/N^2\\
  =\sum_{n\geq 1}\sum_{k_1+\dots+k_m=n}\frac{(-\bu)^\bk}{\bk!}\sum_{g=0}^{\lfloor(n-1)/2\rfloor}
  \sum_{H=(V,E,F)\in\mathcal{H}(g,\bk)}\frac{1}{N^g}\prod_{r=1}^s\la_r^{\#V_r}\prod_{f\in F}\tr_N \mC^{f,\,t}
\\=\sum_{g\geq 0}\frac{1}{N^g}\sum_{k_1+\dots+k_m\geq 1}\frac{(-\bu)^\bk}{\bk!}c(\bk,g,\mC).
\end{multline*}

\end{proof}

%\begin{example}[Example \ref{Ex1}, continued]\label{Ex2}
%  The variance $\rm{Var}  (\tr(q_j))$ is the coefficient
% at $u_j^2$ in \eqref{log Z}. A calculation shows that
%if
%$\mW_1\in\WW{\Sigma_1,p_1},\mW_2\in\WW{\Sigma_2,p_2},\mW_3\in\WW{\Sigma_3,p_3}$
%are independent, then
%  $$
%  \rm{Var}\left(\tr(\mW_1\mW_2)\right)=p_1^2p_2^2\left(\frac{1}{p_1}+\frac{1}{p_2}\right)\tr((\Sigma_1\Sigma_2)^2)+p_1p_2\tr^2(\Sigma_1\Sigma_2);
%  $$
%  \begin{multline*}
%    \rm{Var}\left(\tr(\mW_1\mW_2\mW_3)\right)=
%  p_1^2p_2^2p_3^2\left(\frac{1}{p_1}+\frac{1}{p_2}+\frac{1}{p_3}\right)\tr((\Sigma_1\Sigma_2\Sigma_3)^2)\\
%  +p_1p_2p_3(p_1+p_2+p_3)\tr^2(\Sigma_1\Sigma_2\Sigma_3)+p_1p_2p_3\tr^2(\Sigma_1\Sigma_2\Sigma_3).
%  \end{multline*}
%\end{example}

%\section{Asymptotic normality}\label{Sect CLT}

\section{Proof of Theorem \ref{T3}}\label{Sect CLT}
In order for the method of cumulants (Lemma \ref{Marcinkiewicz}) to work,
we need to reduce the proof to the case of real-valued random variables.
 A sufficient condition for $\tr(q(\mW_1,\dots,\mW_s))$ to be real
is that $q$ is self-adjoint with
respect to the involution $*$ on the space  $\mathbb{C}\langle\bx_1,\dots,\bx_s\rangle$ of complex
non-commutative polynomials, which is defined by
$\bx_1^*=\bx_1,\dots,\bx_s^*=\bx_s$, $a^*=\bar{a}$,
and for $P,Q\in \CC\langle\bx_1,\dots,\bx_s\rangle$, $(PQ)^*=Q^*P^*$.

Suppose $\{q_j: j=1,\dots,d\}$ are the polynomials from the statement of the theorem. If $q$ is one of them,
we can write it as a linear combination of two self-adjoint polynomials with real traces,
$q=\frac12(q+q^*)+\frac12(q-q^*)=\frac12(q+q^*)+\frac1{2i}i(q-q^*)$. So in order to prove the
asymptotic normality of a $\mathbb{C}^d$-valued sequence $(\tr(q_1),\dots,\tr(q_d))$, it suffices  to prove
the asymptotic normality of the $\mathbb{R}^{2d}$-valued
sequence
$$(\tr(q_1+q_1^*),\dots,\tr(q_d+q_d^*),i\tr(q_1-q_1^*),\dots,i\tr(q_d-q_d^*)).$$
To prove the
asymptotic normality of a centered $\mathbb{R}^{2d}$-valued sequence, by Lemma \ref{Marcinkiewicz}
it suffices to show that the variances converge,
and that all cumulants of order higher than $2$ converge to $0$.

By multi-linearity of the cumulants, the joint cumulants $c_\bk$ of the above real random variables,
are the linear combinations of the cumulants
of order $\bk$ of traces of monomials in $\mW_1,\dots,\mW_s$. Therefore, it suffices to show that  all high enough
 cumulants of
(now possibly complex) traces of monomials in $\mW_1,\dots,\mW_s$ vanish in the limit,
and that the second cumulants converge.

\begin{proof}[Conclusion of proof of Theorem \ref{T3}] After the above reduction, it suffices to consider the case,
when
$q_1,\dots,q_d$ are monomials
in $\mW_1^{(N)},\dots,\mW_s^{(N)}$. Denoting $|\bk|=k_1+\dots+k_m$, by Theorem \ref{T2} we have
\begin{equation}
  \label{lim cum}
  c_\bk(\bZ_N)=\sum_{g\geq 0}N^{2-|\bk|-2g}
\sum_{(V,E,F)\in\mathcal{H}(g,\bk)}\prod_{r=1}^s \left(\frac{p_r(N)}{N}\right)^{\# V_r} \prod_{f\in F}\tr_N(\mC^{f,\,t}(N)).
\end{equation}
Thus, $\lim_{N\to\infty}c_\bk(\bX_N)=0$ if $|\bk|\geq 3$.

If  $|\bk|=2$ then  since $p_r(N)/{N}\to \la_r$ and the sum in \eqref{lim cum} is over the finite set of $g\geq 0$,
  we can pass to the limit as $N\to\infty$. We get
\begin{equation}
  \label{lim Var}
  \lim_{N\to\infty}c_\bk(\bZ_N)=\sum_{(V,E,F)\in\mathcal{H}(0,\bk)}\prod_{r=1}^s \la_r^{\# V_r} \prod_{f\in F}m^f,
\end{equation}
where for $f=\{j_1,\dots,j_s\}$, $$m^{f}=\lim_{N\to\infty}\frac1N\tr(\mC_{t(j_1)}(N)\dots\mC_{t(j_s)}(N)).$$
%In the general case, $q_j=\sum_w \alpha_{w,j} \prod_{r\in w}\mW_{t(r)}^{(N)}$, where the sum is taken over a
%finite number of (ordered) finite sequences $w$ of numbers from $\{1,\dots,s\}$.
%This gives
%$$
%\mbox{cov}(\tr (q_i),\tr (q_j))=\sum_{u,w} \alpha_{u,i}\alpha_{w,j}\, \mbox{cov}
%\left(\tr\prod_{r\in u}\mW_{t(r)}^{(N)},\tr\prod_{r\in w}\mW_{t(r)}^{(N)}\right),
%$$
%so the limit $\lim_{N\to\infty}\mbox{cov}(\tr (q_i),\tr (q_j))$ exists for $1\leq i,j\leq d$.
%Since $\overline{\tr(q_j)}=\sum_w \overline{\alpha}_{w,j} \tr\prod_{r\in w^*}\mW_{t(r)}^{(N)}$,
%where $w^*$ is a sequence $w$ in reverse order, so the limit $\lim_{N\to\infty}\mbox{cov}(\tr (q_i),\overline{\tr (q_j)})$
%also exists for $1\leq i,j\leq d$.
\end{proof}

We conclude with deriving a more explicit normalization for a special case. We need the following notation.
Fix $\sigma_1=(1,2,\dots,n)\in\Sn_n$, $\sigma_2=(1,\dots,n)(n+1,\dots,2n)\in\Sn_{2n}$,
, $\sigma_3=(n,n-1,\dots,1)(n+1,\dots,2n)\in\Sn_{2n}$ and  $t:\{1,\dots,n\}\to\{1,\dots,s\}$.
  As previously, let $\Sn_n^\circ(t)\subset\Sn_n(t)$ denote the set of color-preserving permutations $\alpha$ such that
  $\#\sC(\alpha)+\#\sC(\sigma_1\alpha)=n+1$.
  For
  $$\tilde{t}(j)=\begin{cases}t(j)& j\leq n,
  \\
t(j-n)& j>n , \end{cases}
  $$
  let $\Sn_{2n}^*(t)\subset\Sn_{2n}(\tilde{t})$ denote the set of
   permutations $\alpha$ which do not preserve $\{1,\dots,n\}$ and such that
  $\#\sC(\alpha)+\#\sC(\sigma_2\alpha)=2n$. Similarly,
   let $\Sn_{2n}^{**}(t)\subset\Sn_{2n}(\tilde{t})$ consist of
   permutations $\alpha$ which do not preserve $\{1,\dots,n\}$ and such that
  $\#\sC(\alpha)+\#\sC({\sigma}_3\alpha)=2n$.

\begin{corollary}\label{Col CLT}
 Suppose
$\mW_1^{(N)},\dots,\mW_s^{(N)}\in\WW{\mC_N/N,\lceil \la N\rceil }$
are independent, $\la>0$, and there are numbers $m_k$ such that
 $$\lim_{N\to\infty}\tr(\mC_N^{k})-N m_k=0$$ for all $k$.
 Then, as $N\to\infty$,
\begin{equation}\label{VAR}
 \tr(\mW_{t(1)}^{(N)}\dots \mW_{t(n)}^{(N)})-
N\sum_{\alpha\in\Sn_n^\circ(t)} \la^{\# \sC(\alpha)} \prod_{c\in
\sC(\sigma_1\alpha)}m_{\# c}\xrightarrow{\mathcal{D}} Z,
\end{equation}
where $Z=X+iY$, and $X, Y$ are real centered jointly
Gaussian random variables
with the covariance matrix given by:
$$
E(X^2)=\frac12\sum_{\alpha\in\Sn_{2n}^{**}(t)} \la^{\# \sC(\alpha)} \prod_{c\in
\sC({\sigma_3}\alpha)}m_{\# c}+\frac12\sum_{\alpha\in\Sn_{2n}^*(t)} \la^{\# \sC(\alpha)} \Re\left(\prod_{c\in
\sC(\sigma_2\alpha)}m_{\# c}\right),
$$
$$E(Y^2)=\frac12\sum_{\alpha\in\Sn_{2n}^{**}(t)} \la^{\# \sC(\alpha)} \prod_{c\in
\sC({\sigma}_3\alpha)}m_{\# c}-  \frac12\sum_{\alpha\in\Sn_{2n}^*(t)} \la^{\# \sC(\alpha)}\Re\left( \prod_{c\in
\sC(\sigma_2\alpha)}m_{\# c}\right),$$
$$
E(XY)=\frac12 \sum_{\alpha\in\Sn_{2n}^*(t)} \la^{\# \sC(\alpha)} \Im \left(\prod_{c\in
\sC(\sigma_2\alpha)}m_{\# c}\right).
$$
\end{corollary}
\begin{proof} Theorem \ref{T3} implies asymptotic normality, and
Corollary \ref{Col AM} gives the centering, so all that remains is to determine the asymptotic variance.

 We note that if $\alpha\in\Sn_{2n}^*(t)$ then $\langle \alpha,\sigma_2\rangle$ acts transitively on $\{1,\dots,2n\}$,
$\alpha$ preserves $\tilde{t}$, and the hypermap $H=(\sigma_2,\alpha)$ has genus zero.
From  \eqref{lim Var}   with $m=1$ and $q_1(\mW_1\dots\mW_s)=\mW_{t(1)}^{(N)}\dots \mW_{t(n)}^{(N)}$ we see that
\begin{multline*}
E(Z^2)= \lim_{N\to\infty}\mbox{cov}\left(\tr(\mW_{t(1)}^{(N)}\dots
\mW_{t(n)}^{(N)}),\tr(\mW_{t(1)}^{(N)}\dots
\mW_{t(n)}^{(N)})\right)\\=
\sum_{(V,E,F)\in\mathcal{H}(0,2)}\prod_{r=1}^s \la^{\# \sC_r(\alpha)} \prod_{f\in F}m^f
=
\sum_{\alpha\in\Sn_{2n}^*(t)} \la^{\# \sC(\alpha)} \prod_{c\in
\sC(\sigma_2\alpha)}m_{\# c}.
\end{multline*}
Similarly,
\begin{multline*}
  E(|Z|^2)=\lim_{N\to\infty}\mbox{cov}\left(\tr(\mW_{t(1)}^{(N)}\dots
\mW_{t(n)}^{(N)}), \tr(\mW_{t(n)}^{(N)}\dots
\mW_{t(1)}^{(N)}) \right)\\=\sum_{\alpha\in\Sn_{2n}^{**}(t)} \la^{\# \sC(\alpha)} \prod_{c\in
\sC({\sigma}_3\alpha)}m_{\# c}.
\end{multline*}
This ends the proof, as $E(X^2)=E(|Z|^2)/2+\Re E(Z^2)/2$, $E(Y^2)=E(|Z|^2)/2-\Re E(Z^2)/2$, and
$E(XY)=\Im E(Z^2)/2$.
\end{proof}
%\begin{remark} The proof implies that
%$$\sum_{\alpha\in\Sn_{2n}^{**}(t)} \la^{\# \sC(\alpha)} \prod_{c\in
%\sC({\sigma}_3\alpha)}m_{\# c}\geq \Re\sum_{\alpha\in\Sn_{2n}^*(t)} \la^{\# \sC(\alpha)} \prod_{c\in
%\sC(\sigma_2\alpha)}m_{\# c}$$
%for any $\la>0$ and any sequence $(m_k)$ of moments of a probability measure on $[0,\infty)$.
%%\fbox{True or False?} If $\alpha$ does not preserve $\{1,\dots,n\}$
%%then $g(\sigma_2,\alpha)\geq g(\sigma_3,\alpha)$.
%\end{remark}
\begin{example}
  \label{Ex CLT} Here we illustrate the fact that under additional assumptions one can calculate explicit
 normalization for Theorem \ref{T3}.
From Examples \ref{Ex1} and \ref{Ex2}, we see that if
$\mW_1,\mW_2,\mW_3\in \WW{\mC_N/N,\lceil \la N\rceil }$
are independent,  $\lim_{N\to\infty}\tr_N(\mC_N^{k})= m_k$ for all $k$, and
$\lim_{N\to\infty}\tr(\mC_N^{3})-N m_3=b$, then
$$
\tr(\mW_1\mW_2\mW_3)- N \la^3 m_3 \xrightarrow{\mathcal{D}} X+iY  \mbox{ as $N\to\infty$},
$$
where $X,Y$ are independent real Gaussian random variables;
$X$ has mean $\la^3 b$ and variance $3\la^5m_6+3\la^4(m_2m_4+m_3^2)/2+\la^3m_2^3/2$;
$Y$ has mean zero and variance $3\la^4(m_2m_4-m_3^2)/2+\la^3m_2^3/2$.
In particular, if $\mC_N=I$,  then $X,Y$ have mean zero and variances $\la^3(3\la^2+3\la+1/2)$
and $\la^3/2$ respectively.
\end{example}

\subsection*{Acknowledgement} The author thanks J. Mingo for several helpful discussions and references.
He also thanks A. Guionnet for bringing to his attention Wishart matrices, and to H. Massam for pointing
out the importance of Wishart matrices with general covariances. Finally, he thanks the anonymous
referee for the exceptionally detailed report.

%BibTeX
%\bibliographystyle{alpha}
%\bibliographystyle{apalike}
\bibliographystyle{newapa} %name (year)
\bibliography{matrix-models}
\end{document}